\documentclass{elsart}

\usepackage{amsfonts}
\usepackage{graphicx}
\usepackage{amsmath,natbib,subfigure}

\def\R{\mathbb R}

\def\epsilon{\varepsilon}
\def\ds{\displaystyle}

\newcommand{\be}{\begin{equation}}
\newcommand{\ee}{\end{equation}}

\begin{document}

\begin{frontmatter}

\title{Does reaction-diffusion support the duality of fragmentation effect? }

\author[aut1]{Lionel Roques\corauthref{cor}},
\ead{lionel.roques@avignon.inra.fr}
\author[aut2]{and M.D. Chekroun}
\corauth[cor]{Corresponding author. Fax: +33 432722182.}
\address[aut1]{UR 546 Biostatistique et Processus Spatiaux, INRA, F-84000 Avignon, France}
\address[aut2]{Environmental Research and Teaching Institute, \'Ecole
Normale Sup\'erieure,\\ 24 rue Lhomond, 75231 Paris Cedex 05,
France}

\begin{abstract}
There is a gap between single-species model predictions, and
empirical studies, regarding the effect of habitat fragmentation
per se,  i.e., a process involving the breaking apart of habitat
without loss of habitat. Empirical works indicate that
fragmentation can have positive as well as negative effects,
whereas, traditionally, single-species models predict a negative
effect of fragmentation. Within the class of reaction-diffusion
models, studies almost unanimously predict such a detrimental
effect. In this paper, considering a single-species
reaction-diffusion model with a removal -- or similarly harvesting
-- term, in two dimensions, we find both positive and negative
effects of fragmentation of the reserves, i.e. the protected
regions where no removal occurs. Fragmented reserves lead to
higher population sizes for time-constant removal terms. On the
other hand, when the removal term is proportional to the
population density, higher population sizes are obtained on
aggregated reserves, but maximum yields are attained on fragmented
configurations, and for intermediate harvesting intensities.
\end{abstract}


\end{frontmatter}


\section{Introduction}\label{sec_intro}

The analysis of the effects of environmental fragmentation and
variability
 on population densities and biodiversity has stimulated
the development of many spatially-explicit population models. In
the modeling literature, positive effects of environmental
variability have been recorded \citep[see e.g.][]{bolker,bhr1}. On
the other hand \cite{fahrig2003}, in a thorough bibliography
analysis, pointed out that most single-species modeling approaches
lead to comparable conclusions regarding the
 detrimental effects of fragmentation per se,  i.e.,
a process involving the breaking apart of habitat without loss of
habitat. She noted that unlike the effects of habitat loss
\citep[see][for a discussion of the consequences of fragmentation
with habitat loss]{shm}, and in contrast to current theory,
empirical studies suggest that the effects of fragmentation  per
se are at least as likely positive as negative. The aim of this
note is to make steps towards a reconciliation between the theory
and empirical works on the effects of fragmentation, within the
framework of reaction-diffusion models.

Reaction-diffusion models (hereafter RD models), although they
sometimes bear on simplistic assumptions such as infinite velocity
assumption and completely random motion of animals
\citep{Holmes93}, are not in disagreement with certain dispersal
properties of populations observed in natural as well as
experimental ecological systems, at least qualitatively
\citep[see][]{sk,tur,mu,oku}. Furthermore, these models often
provide a good framework for rigorous investigation of theoretical
questions and derivation of qualitative as well as numerical
results on population dynamics. In that respect, the effects of
environmental fragmentation have been addressed in many
theoretical studies based on such models over the last decades.

Within the class of RD models with heterogenous coefficients,
numerous recent works have emphasized the detrimental effect of
environmental fragmentation  per se on species persistence and
spreading, in agreement with the other theoretical tendencies
noted
 by Fahrig. In all these RD models, the population growth rate
function at a location $\mathbf{x}$, $r(\mathbf{x})$, was not
constant, taking higher values in favorable regions than in
unfavorable ones. Depending on the spatial arrangements of these
regions, the modeled populations were shown to tend to extinction
or survive, and to disperse at different speeds. In the particular
case of 1-dimensional binary environments (i.e. for $r$ taking
values in a set constituted of two values), \citet{cc0,ccL} and
\cite{sk} have established that concentrating all the habitat in a
single patch improved persistence. \cite{bhr1} generalized these
analytical results to the $N$-dimensional case, with more general
growth rate functions. More recently \cite{rs1} carried out more
precise results regarding the negative correlation between
persistence and fragmentation. Habitat fragmentation has also been
shown, first numerically \citep{kkts}, and then analytically
\citep{ehr} to negatively affect population spreading speed.

Nevertheless, to our knowledge, there do not exist detailed
demonstration of positive effects related to fragmentation derived
from RD models in two space-dimensions, although in a related
work,  \cite{neub} shows positive effects of reserve fragmentation
from RD models with harvesting terms, in the sense that fragmented
reserves sometimes maximize the yield. However, in his work, which
is carried out in a one-dimensional space, no explicit
conservation of the area of the reserve is assumed.

In this paper, we study  single-species RD models in two-space
dimensions with a spatially-homogeneous growth term, and spatially
heterogeneous removal terms $-Y$ which can be, for instance,
interpreted as harvesting terms. The regions where $Y\equiv 0$
hence correspond to protected regions or similarly reserves. In
our models, $Y=Y(\mathbf{x},u)$ depends on the location
$\mathbf{x}$, and can depend on the population density $u$ at this
location. Typically, when $Y$ does not depend on the population
density, it corresponds to a \emph{constant-yield harvesting}
strategy. In this case, a constant number of individuals are
removed per unit of time. This is the case when a quota is set on
the harvesters \citep{rob1,rob2,steph}. Even in the absence of
such imposed quotas, harvesters often increase their effort to
maintain a constant yield. A good example is provided by the
high-trophic level fishes catches data in the North Atlantic
reviewed by \cite{christ}, which describe a decline of the biomass
of one half from 1950 to 1990, while the catch remains the same.
The function $Y$ can also be taken to be proportional to the
density $u$, corresponding to a \emph{proportional harvesting}
strategy. Then, locally, a constant proportion of the population
is removed per unit of time, corresponding to a constant effort of
the harvesters. These two harvesting strategies have been
investigated by the authors  \citep{rc1},  in inhomogeneous
environments. It was shown, through analytical and numerical
analysis, that aggregated habitat configurations gave better
chances of population persistence, respecting the tendency found
in the modeling literature; see also \cite{shi2} for other
mathematical results on these models.

We conduce here a different analysis. Not only we investigate the
effects of the spatial arrangement of the harvesting term rather
than those of the growth function, but we also focus on other
quantities than simple persistence. Considering protected regions
with a fixed total area but with gradually fragmented shapes, we
analyze the intertwined effects of fragmentation and harvesting
intensity on both the population size and the quantity of
harvested individuals. To do that we use the stochastic model of
landscape generation of \cite{rs1}, that is in complete agreement
with the concept of fragmentation per se, and we show that
fragmentation of protected regions can in fact be beneficial to
the modeled population and to the harvesters. These results
demonstrate that two-dimensional RD models with harvesting terms
can support a dual effect of fragmentation per se, positive as
well as negative.

Typically, it is shown, for instance, that fragmented reserves
lead to higher population sizes for time-constant removal terms.
On the other hand, when the removal term is proportional to the
population density, higher population sizes are obtained on
aggregated reserves, but maximum yields are attained on fragmented
configurations, and for intermediate harvesting intensities.

\section{Materials and methods}\label{sec_MM}

\subsection{The model}\label{sec_model}

The idea of modeling population dynamics with reaction-diffusion models has begun
to develop at the beginning of the 20th century, with random walk
theories of organisms, introduced by \cite{pm}. Then,  \cite{fi} and \cite{kpp}
used a reaction-diffusion equation with homogeneous coefficients as
a model for population genetics. Later, \cite{ske} examined
this type of model, and he succeeded to propose quantitative
explanations of observations for the spread of muskrats throughout
Europe at the beginning of 20th Century. Since then, these models
have been widely used to explain spatial propagation or spreading of
biological species  \citep[bacteria, epidemiological agents, insects,
fishes, mammal, plants, etc., see the books][for review]{sk,tur,mu,oku}.

Ignoring age or stage structures as well as delay mechanisms or
Allee effects, the classical Fisher-Kolmogorov model, in two
space-dimensions, can be written as follows: \par\nobreak\noindent
\begin{equation}\label{eqkpp}
\frac{\partial u}{\partial t}-D \nabla^2 u=r  u(1- u/K), \ t>0, \ \mathbf{x}\in\Omega\subset\R^2,
\end{equation}
where $u=u(t,\mathbf{x})$ corresponds to the population density at
time $t$ and position $\mathbf{x}=(x_1,x_2)$. The left-hand side of
(\ref{eqkpp}) corresponds to the diffusion equation, and simply
describes the redistribution of organisms following uncorrelated
random walks where $\nabla^2$ stands for the spatial dispersion
operator $\ds{\nabla^2 u=\frac{\partial^2 u}{\partial
x_1^2}+\frac{\partial^2 u}{\partial x_2^2}}$. The diffusion
coefficient  $D$ measures the individuals rate of movement, $r>0$ is
the intrinsic growth rate of the population and $K>0$ corresponds to
the environment carrying capacity.

The domain $\Omega$ is considered bounded, and we assume
reflecting boundary conditions: \par\nobreak\noindent
$$\frac{\partial u}{\partial n}(t,\mathbf{x})=0, \hbox{ for }\mathbf{x}\in \partial
\Omega,$$where $\partial \Omega$ is the domain's boundary and
$n=n(\mathbf{x})$ corresponds to the outward normal to this
boundary. Thus, some part of the boundary  can be considered as an
absolute barrier that the individuals do not cross, like coasts, and
other parts of the boundary can be seen as regions where as much
individuals exit the domain as individuals enter the domain.

At this stage no environmental fragmentation is present in this
class of models. To introduce it, we adopt a perturbative approach
which consists in subtracting a spatially-dependent
 term to the right-hand side of equation (\ref{eqkpp}): \par\nobreak\noindent
\begin{equation}
\frac{\partial u}{\partial t}-D \nabla^2 u=r  u(1-
u/K)-Y(\mathbf{x},u), \ t>0, \ \mathbf{x}\in\Omega\subset\R^2.
\label{eq_base}
\end{equation}

It could seem artificial to keep constant the parameters $r,D$ and
$ K$ whereas the removal term $Y$ is spatially-dependent. The
ecological interpretation of such a framework is no more than the
consideration of heterogeneously distributed harvesting in
homogeneous media, a classical set-up in fisheries for instance.
In fact, even in the case of a population living in a
heterogeneous environment, the objective of such assumptions is to
separate the effects in order to facilitate the ecological
interpretations; and keeping constant the biological parameters
can be thought as to consider the effects of spatially-dependent
perturbations on an averaged model of growth (\ref{eqkpp}), where
the parameters are averaged in space and time. Advanced
mathematical theory of averaging for partial differential
equations can then be used to support this strategy and to make
robust conclusions derived by such approach, at least in the case
of small amplitude of oscillating parameters
$r(t,\mathbf{x}),D(t,\mathbf{x})$ and $K(t,\mathbf{x})$ \citep[see
e.g.][]{hale, cr1} .

In the forthcoming computations, we assume that harvesting starts
on a previously not harvested population, which has reached its
stable positive steady state. In other words, the environment is
assumed to have reached its carrying capacity at $t=0$, the time
at which harvesting is started, with $u(0,\mathbf{x})=K$, where
$u$ is the solution of (\ref{eq_base}). This assumption is a
natural one for studying the destabilizing effects due to removal
terms. Our model (\ref{eq_base}) being introduced and discussed,
the analysis of environmental effects brought by the removal term
$Y$ is addressed according two harvesting strategies that we
present now.

\subsection{Harvesting strategies \label{sec_strateg}}

The first type corresponds to  a ``quasi"-constant-yield
harvesting, with a removal term
\par\nobreak\noindent \begin{equation} Y(\mathbf{x},u)=\delta \cdot
\chi(\mathbf{x}) \rho_{\varepsilon}(u). \label{const_yield}
\end{equation} This model is close to but different from the
so-called threshold harvesting model found in ordinary
differential equations (ODEs) models of harvesting
\citep[see][]{rc1}.

In  (\ref{const_yield}), $\chi(\mathbf{x})$ is a function
taking the value $1$ if $\mathbf{x}$ belongs to a harvested
region, and the value $0$ if $\mathbf{x}$ belongs to a protected
region. It constitutes therefore what we call the harvesting
field, with $\delta$ a positive constant which corresponds to the
harvesting intensity in this field. The parameter $\delta$ can
also be interpreted as a quota. The last term
$\rho_{\varepsilon}(u)$ is a density-dependent threshold function:
\par\nobreak\noindent
$$\rho_{\varepsilon}(s)=0 \, \hbox{ if } s\le 0, \
\rho_{\varepsilon}(s)=s/\varepsilon \, \hbox{ if }
0<s<\varepsilon, \hbox{ and } \rho_{\varepsilon}(s)=1 \ \hbox{ if
} s\geq \varepsilon,$$ where $\varepsilon$ is a small threshold
below which harvesting is progressively withdrawn. With such a
harvesting function, at each location $\mathbf{x}$, the yield is
constant in time whenever $u(t,\mathbf{x}) \geq \varepsilon$. Note
that from a mathematical point of view, the function
$\rho_{\epsilon}$ ensures the nonnegativity of the solutions of
 (\ref{eq_base}) \citep{rc1}. Considering constant-yield
harvesting functions without this threshold value would be
unrealistic since it would lead to harvest on zero-populations.

The second type of strategy corresponds to a more standard
proportional harvesting situation \citep{neub}, in a spatial
context, where \par\nobreak\noindent \be Y(\mathbf{x},u)=E \cdot
\chi(\mathbf{x}) u. \label{prop_harv} \ee The function $\chi$ is
defined as above, and the term $E \cdot \chi(\mathbf{x})$ can now
be interpreted as a harvesting effort at the location
$\mathbf{x}$. The instantaneous yield at a point $\mathbf{x}$ is
then proportional to this effort and to the local population
density.

\subsection{The model of fragmentation}\label{sec_frag}

There exist several ways of obtaining hypothetical landscape
distributions, see  e.g.   \citet{gard} and \citet{keitt} for
neutral landscape models, and \citet{mandel} for measures of
fragmentation based on fractal dimension. The model retained here
is the one developed by \citet{rs1}, and inspired from statistical
physics. This is a neutral landscape model in the sense that it is
a stochastic model of landscape pattern, and the value --
protected or harvested in the present case -- assigned to a
position in the pattern is a random variable. As key property,
this stochastic model  provides a numerical procedure for
generating several samples of landscapes with breaking apart of
the habitat while keeping constant habitat abundance, along with
an exact control of the later and the type of breaking apart. This
model offers therefore an appropriate framework for assessing the
effect of fragmentation  per se on RD models, in the sense
underlined by \cite{fahrig2003}. We make precise here the main
parameters calibrating the model for our present purpose.

\begin{figure} \centering
\subfigure[$s(\chi)=94$]{%
\includegraphics*[width=4cm]{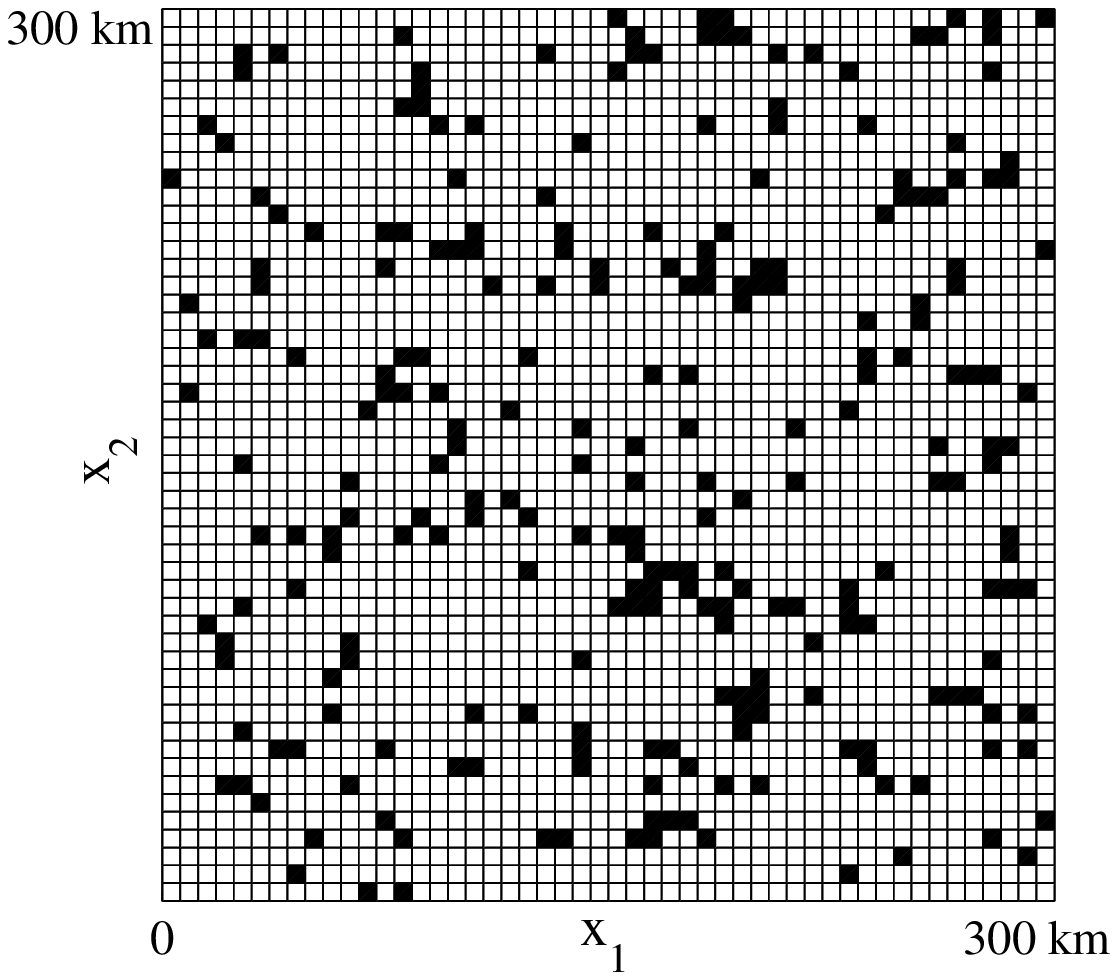}}
\subfigure[$s(\chi)=274$]{%
\includegraphics*[width=4cm]{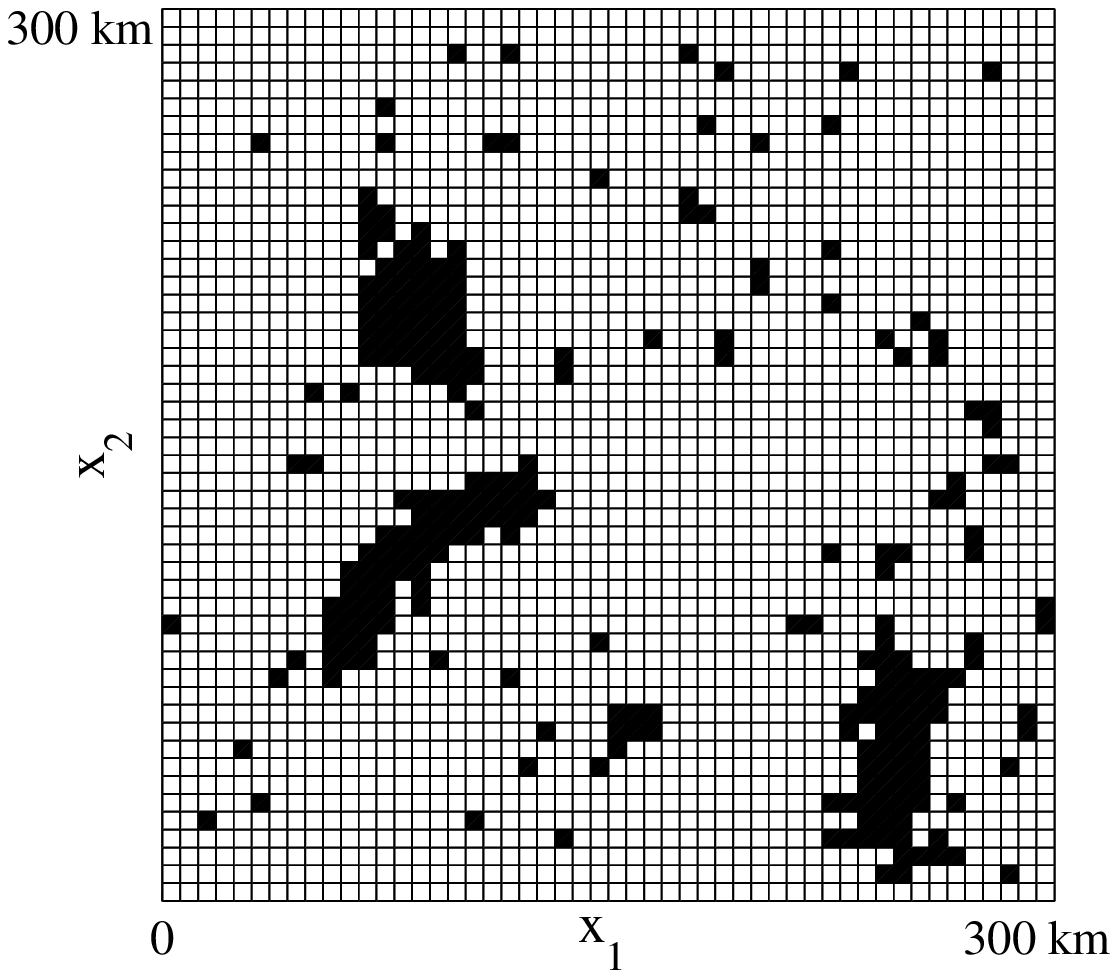}}
\subfigure[$s(\chi)=460$]{%
\includegraphics*[width=4cm]{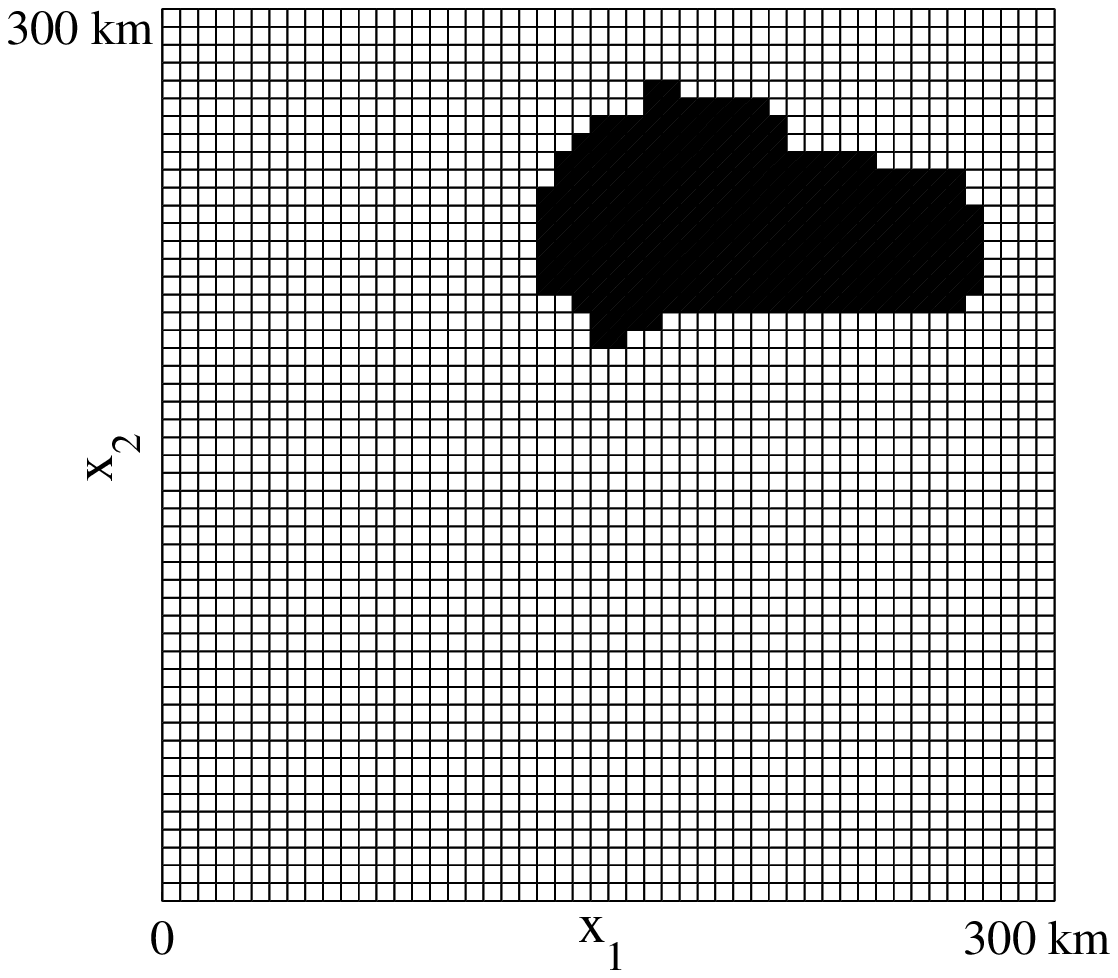}}
\caption{\footnotesize Some samples of harvesting field
configurations.
The black areas correspond to protected regions, where $\chi=0$.}%
\label{fig:frag} 
\end{figure}

The harvested and protected regions are entirely determined by the
harvesting field $\chi(\mathbf{x})$ present in (\ref{eq_base}) via
the removal term $Y$. In order to build gradually fragmented
configurations of these regions, we have discretized the domain
$\Omega$ into $50\times 50$ subcells $C_i$, with in some cells
$\chi(\mathbf{x})=0$ for the unharvested cells, and
$\chi(\mathbf{x})=1$ in the other cells. Based on the stochastic
model of \cite{rs1}, we have built $6000$ samples of such
functions $\chi(\mathbf{x})$, with different degrees of
fragmentation. In all these samples, the protected region occupies
10 $\%$ of the domain $\Omega$. The fragmentation of the protected
region is defined as follows. The lattice made of the cells $C_i$
is equipped with a 4-neighborhood system $V(C_i)$. We set
$s(\chi)=$number of pairs of neighbors $(C_i,C_j)$ such that
$\chi$ takes the value $0$ on $C_i$ and $C_j$. This number
$s(\chi)$ is directly linked to fragmentation: the protected
region is all the more aggregated as $s(\chi)$ is high, and all
the more fragmented as $s(\chi)$ is small (Fig. \ref{fig:frag}).
Therefore,  $s(\chi)$ can be seen as an ``aggregation index'' of
the protected region. On our samples, the aggregation index $s$
varies from 94 to 460. For each aggregation index incremented as
follows $s_k:=94+6 \times (k-1)$, we picked up arbitrarily a
configuration $\chi_k$ with $s(\chi_k)=s_k$. This lead to $62$
harvesting field distributions, with gradually aggregated
configurations of the protected region.

\subsection{Methods}\label{sec_method}

The response to the spatial perturbation terms $Y(\mathbf{x}, u)$,
distributed according to the aggregation index of their underlying
harvesting fields, as described above, is analyzed in terms of
total population size $P(t)$, and annual yield $R(t)$ in the
region $\Omega$. More precisely, we evaluate the time-dependent
quantities \par\nobreak\noindent
$$P(t):=\int_{\Omega}u(t,\mathbf{x})\mbox{d}\mathbf{x},$$
corresponding to the total population at time $t$, and
\par\nobreak\noindent
$$R(t):=\int_{t-1}^{t}\int_{\Omega}
Y(\mathbf{x},u(\tau,\mathbf{x}))\mbox{d}\mathbf{x} \mbox{d}
\tau,$$ corresponding to the annual yield during the year that
precedes $t$.

\textit{Note 1:}  In the case of quasi-constant-yield harvesting,
if $u(\tau,\mathbf{x})$ is greater or equal than $\varepsilon$ for
all $\tau$ living within the temporal window $(t-1,t)$ and
$\mathbf{x}$ in the domain $\Omega$, we simply obtain
\par\nobreak\noindent
$$R(t)=\delta \cdot [\hbox{area of the harvested region}].$$In the
case of proportional harvesting, we have: \par\nobreak\noindent
$$R(t)=E \cdot [\hbox{mean population in the harvested region
during the year }(t-1,t)].$$

We estimate the intertwined effects of fragmentation  per se and
harvesting intensity by plotting $P(t)$ and $R(t)$ against
$\delta$ for the quasi-constant yield strategy, and against $E$
for the proportional harvesting strategy. For every configuration,
the result is represented by a curve which color is attributed in
function of the aggregation index; red corresponds to the more
aggregated configurations, and blue to the more fragmented
configurations. For each fixed harvesting intensity, we computed
the gap between the maximum and minimum population sizes obtained
over the $62$ harvesting field distributions, and we expressed it
in terms of relative loss obtained in the worst configuration
compared to the best one, through the formula: $100 \times
(\hbox{highest population}-\hbox{lowest
population})/(\hbox{highest population})$. Similarly, for the
annual yield, $100 \times (\hbox{highest annual
yield}-\hbox{lowest annual yield})/(\hbox{highest annual yield})$
was computed for each fixed harvesting intensity, where the maxima
and minima are taken over the $62$ harvesting field distributions.

For the numerical setup we consider $\Omega$ to be a square domain
of $300$~km$\times 300$~km. We set $r=1$~year$^{-1}$ and
$K=10^3$~individuals/km$^2$. The diffusion coefficient $D$ varies
between $10$~km$^2$/year (low mobility) and 100~km$^2$/year (high
mobility); see the book of \citet{sk} for some observed values of
$r$ and $D$, for several animal species. The threshold
$\varepsilon$ is set to $10$~individuals/km$^2$. Our results, except in Section \ref{sec_resPt}, are presented at a fixed time $t=5$;
this time has been chosen to fit usual times for observations and responses in anthropic harvesting activities.

The numerical integrations of the RD models were performed using a
second order finite elements method where the solutions
$u(t,\mathbf{x})$ of the model (\ref{eq_base}) have been computed
with the initial condition and harvesting strategies discussed
above. The quantities of interest $P(t)$ and $Q(t)$ are then easily
computable. The numerical results with ecological interpretations
are discussed in the following section.

\section{Results}\label{sec_result}

\subsection{Quasi-constant-yield harvesting
strategy}\label{sec_resconst}

\begin{figure}
\centering
\includegraphics*[width=14cm]{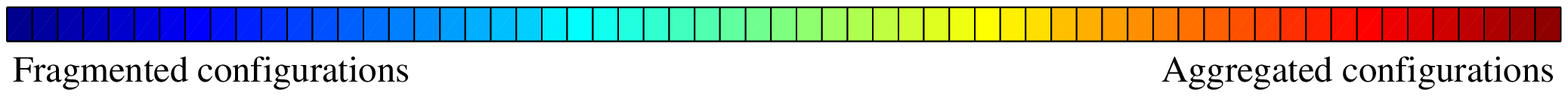}
\subfigure[]{%
\includegraphics*[width=6.5cm]{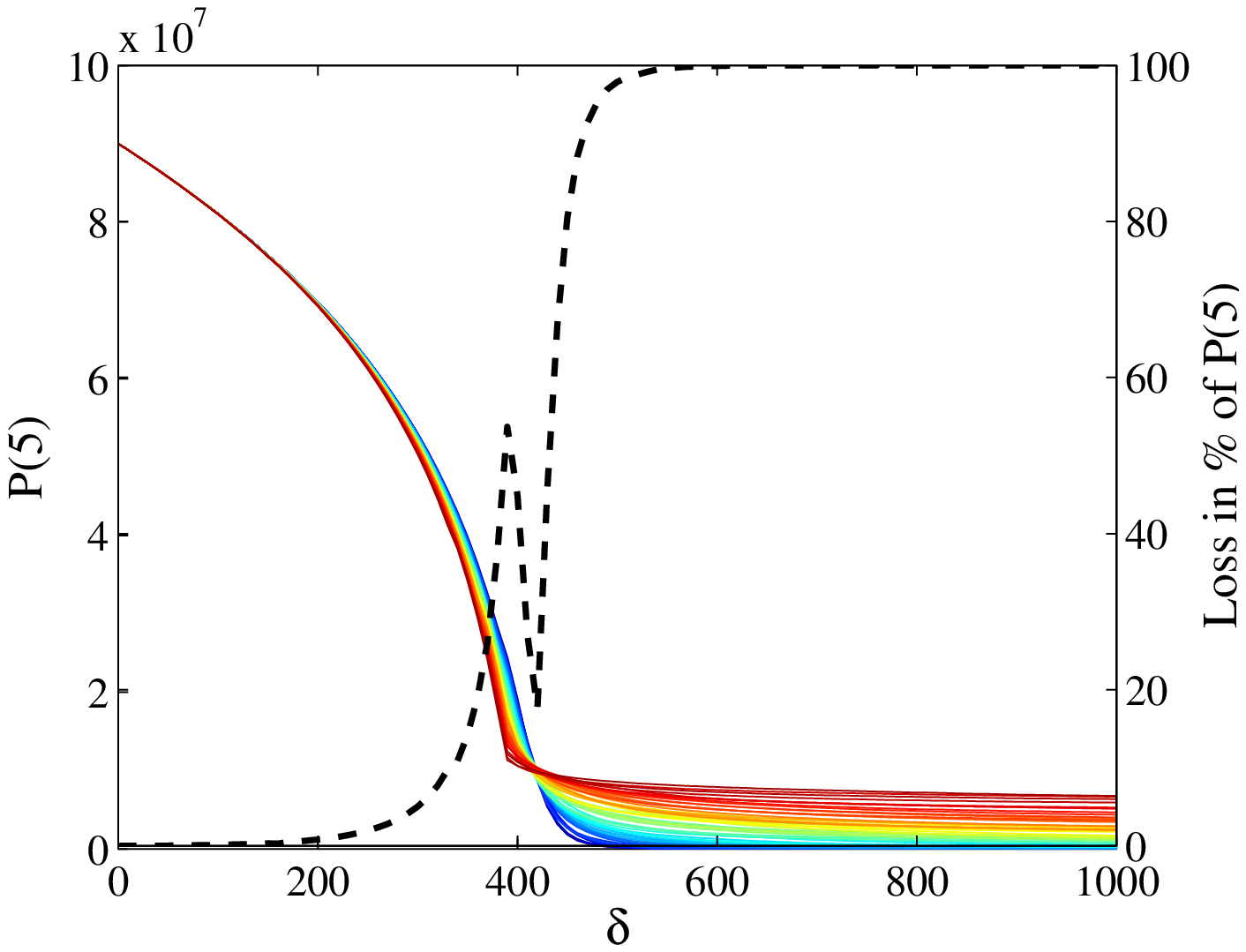}}
\subfigure[]{%
\includegraphics*[width=6.5cm]{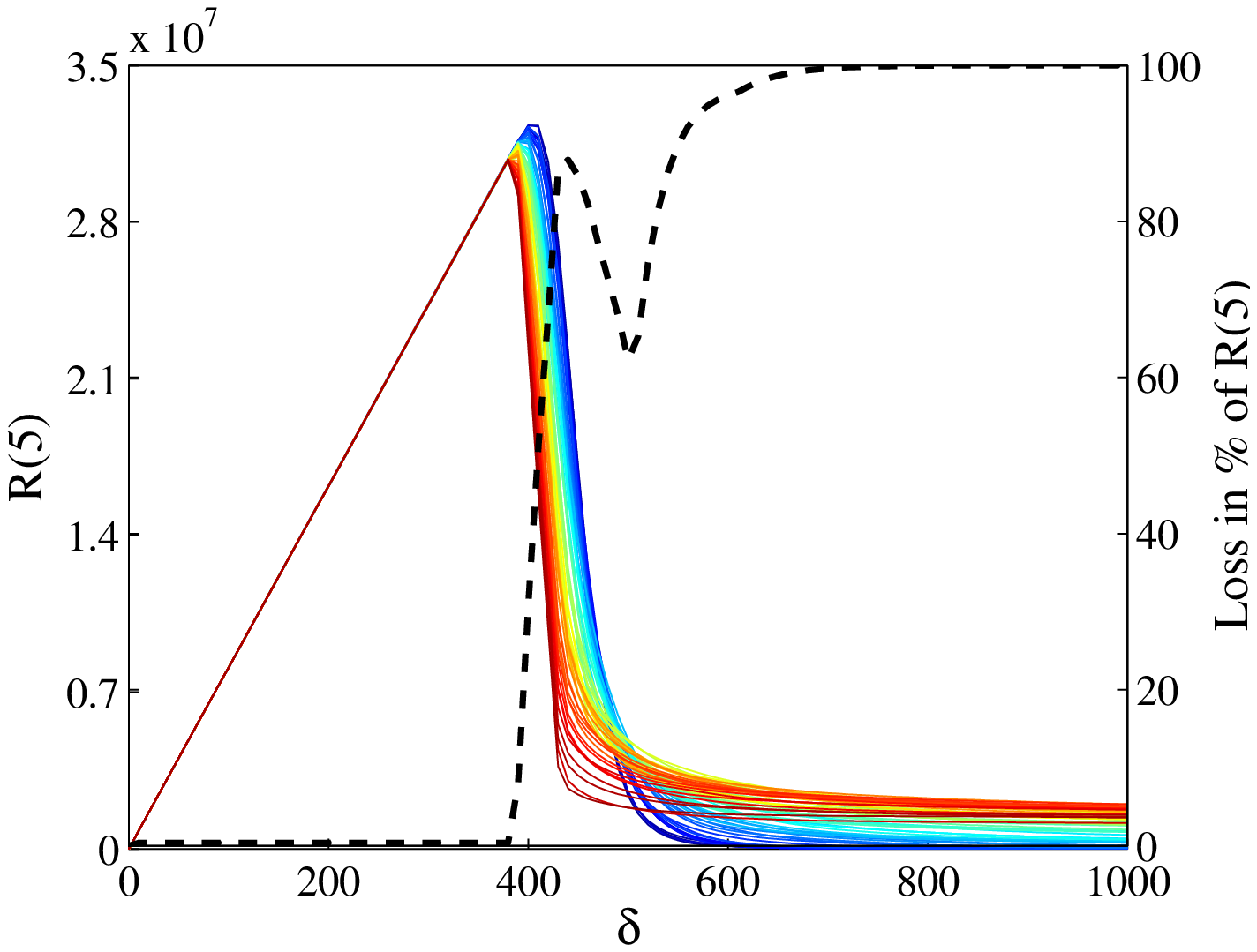}}
\caption{\footnotesize (a): Total population in $\Omega$ after $5$
years, in function of the quota $\delta$. (b): Total yield during
year $4$, $R(5)$, in function of $\delta$. Each curve is associated
with a different configuration of the protected region. Blue curves
correspond to more fragmented configurations, and red curves
correspond to more aggregated configurations. The black dotted lines
indicate, for each $\delta$,  the relative losses obtained in the
worst configurations compared to the best ones.
These computations were carried out for $D=50$.}%
\label{fig:const}
\end{figure}

For every configuration, the higher the quota $\delta$, the smaller
the  population size $P(5)$ (Fig. \ref{fig:const}, a). On the other
hand, a maximum yield $R(5)$ is reached for an intermediate value of
the quota, while small values and large values of $\delta$ both lead
to small yields (Fig. \ref{fig:const}, b).

As it could be expected, for small values of $\delta$, there is no
dependence of the yield with respect to the habitat configuration.
Indeed, in such cases, the population density should be everywhere
above $\varepsilon$, and the harvesting function is therefore
constant in time, equal to $\delta$ in the harvested regions; see
Note 1. However, this linear dependence of $R(5)$ pursues for
higher values of $\delta$ in fragmented configurations. This
indicates that the population density  never falls below
$\varepsilon$ on these configurations. Thus, the more fragmented
configurations first lead to higher yields for intermediate
quotas, with a yield loss in aggregated configurations which
attains $88\%$. Conversely, the more fragmented configurations
lead to lower yields for higher values of $\delta$. Overall, the
maximum yield is attained on the more fragmented configuration.
Similarly, population sizes are higher for fragmented
configurations and low quotas, with losses in population size up
to $54\%$ on aggregated configurations. They then become higher
for aggregated configurations and higher quotas. It is noteworthy
that, whenever the harvesting term is really constant, i.e., in
the region where $R(5)$ is linear, fragmented configurations lead
to higher population sizes.

Such a reversion of the influence of fragmentation on both $P$ and
$R$, for increasing quotas, is not intuitively obvious. Yet, we can
give a reasonable explanation for it. In fragmented configurations,
the mean distance to a protected region --- where population density
is higher --- is reduced compared to more aggregated configurations.
Therefore, when $\delta$ is not too large, at each location in
$\Omega$, the population can be efficiently sustained by the
protected regions, and the density never falls below the threshold
$\varepsilon$, leading to higher yields, compared to more aggregated
configurations. For larger values of $\delta$, the harvested regions
become very hostile, and because of dispersion, populations tend to
extinct, even in the protected regions. In the case of protected
regions with small perimeters, corresponding to aggregated
configurations, dispersion of the individuals into the harvested
regions is reduced.  With such configurations, populations can
therefore sustain higher quotas without risking extinction, leading
to higher values of $P$ and $R$. Note that, whenever
$u(t,\mathbf{x})$ is less than or equal to  $\varepsilon$ everywhere
in the harvested region, quasi-constant-yield harvesting becomes
equivalent to proportional harvesting, with effort
$E=\delta/\varepsilon$.

Comparable qualitative results were found for diffusion
coefficients $D$ ranging from $10$ to $100$ and are thus not
shown.

\subsection{Proportional harvesting strategy}\label{sec_resprop}

When the harvesting function $Y$ is of proportional type
(\ref{prop_harv}), our model reduces to $\frac{\partial
u}{\partial t}-D \nabla^2 u=r u(1-E \cdot \chi(\mathbf{x})/r-
u/K)$. Population persistence for this model has been thoroughly
investigated \citep{ccL,bhr1,rh1,rs1}, both analytically and
numerically.

The specific effects of fragmentation of the protected regions can
be deduced from the numerical study of \cite{rs1}. It shows that, on
aggregated configurations, higher efforts $E$ can be sustained
without risking extinction. The optimal shapes of the protected
regions, in terms of maximum sustainable effort, has even been
obtained in \citep{rh1}. Surprisingly, these shape depend on the area
of the protected region. Indeed, small areas have been proved to
lead to disc-shaped optimal shapes, while high areas lead to
stripe-shaped optimal shapes.

\begin{figure}
\centering
\includegraphics*[width=14cm]{colorbar.eps}
\subfigure[]{%
\includegraphics*[width=6.5cm]{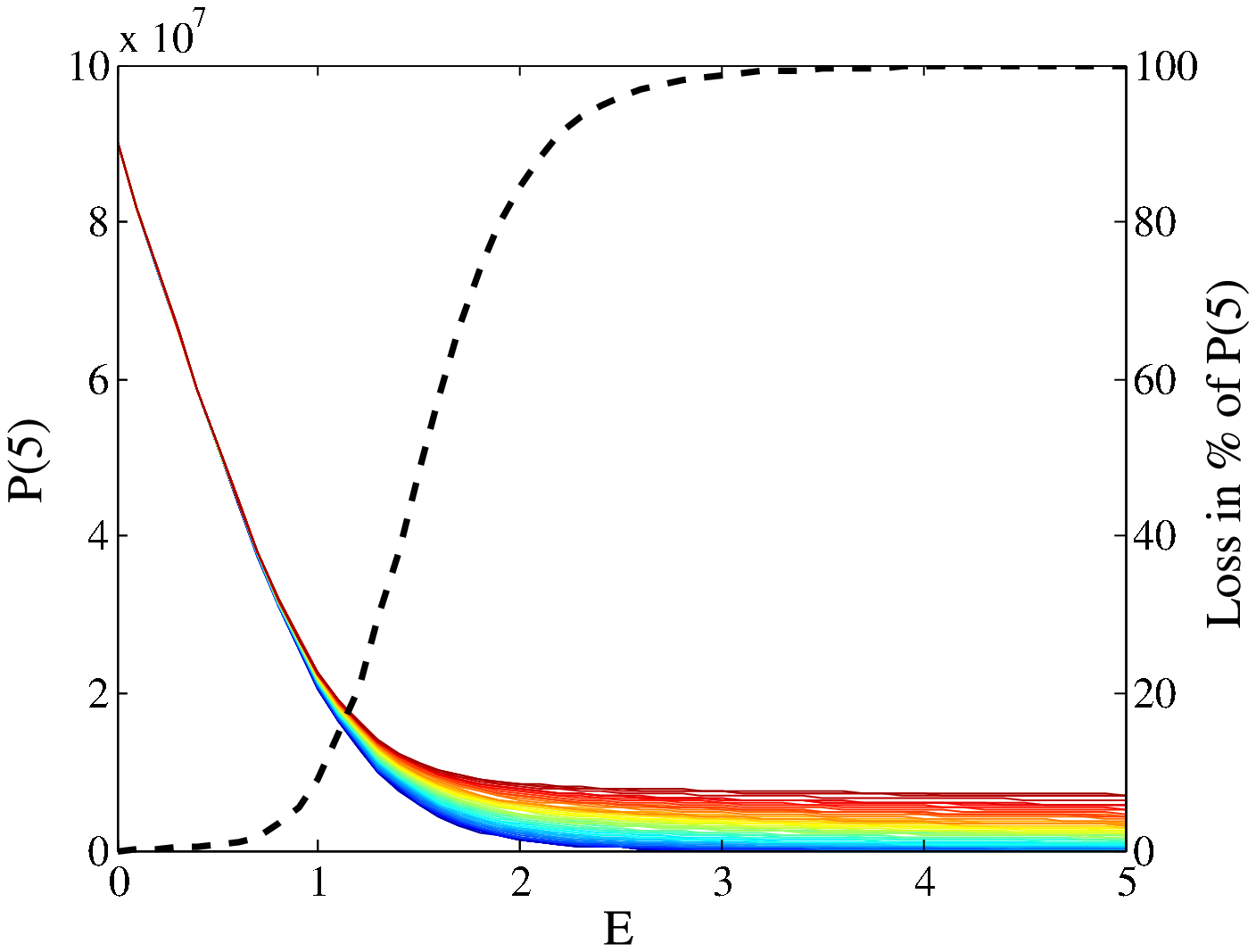}}
\subfigure[]{%
\includegraphics*[width=6.5cm]{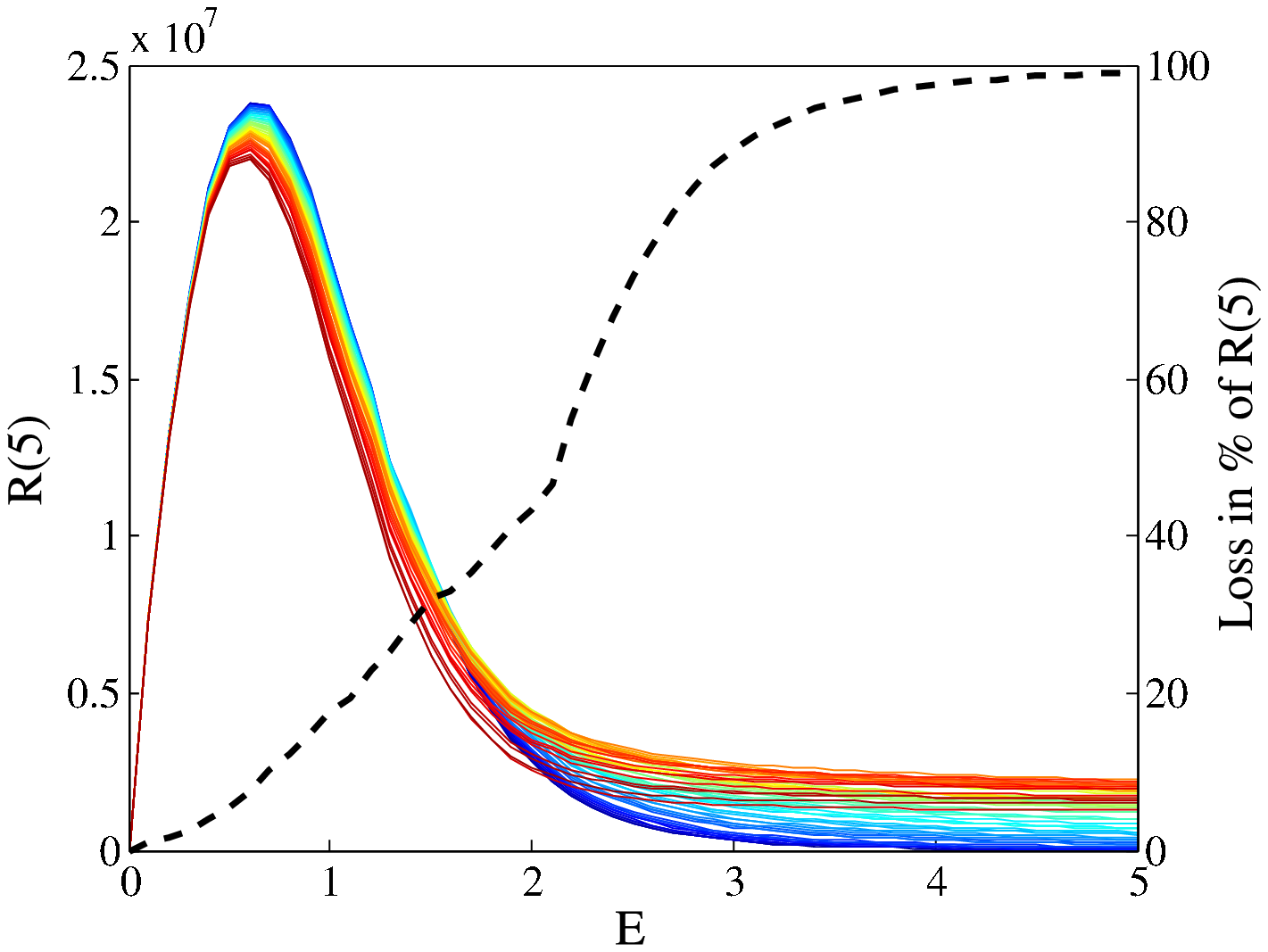}}
\caption{\footnotesize (a): Total population in $\Omega$ after $5$
years, in function of the harvesting effort $E$. (b): Total yield
during year $4$, $R(5)$, in function of $E$. To each curve
corresponds a different configuration of the protected regions; the
blue curves correspond to more fragmented configurations, and the
red curves correspond to more aggregated configurations. The black
dotted lines indicate, for each value of $E$,  the relative losses
obtained in the worst configurations compared to the best
ones. For these computations, we fixed $D=50$.}%
\label{fig:prop}
\end{figure}

Following our approach, we still focus on the quantities $P$ and
$R$, on which the effects of fragmentation have not yet been
investigated. For every configuration, the population size $P(5)$
decreases as the effort $E$ increases (Fig. \ref{fig:prop}, a).
Contrarily to the quasi-constant case, with such a proportional
harvesting strategy, aggregated configurations always lead to
larger populations, whatever the effort. The effect of
fragmentation/aggregation of the protected region on the
population sizes really becomes noticeable when $P(5)$ falls below
one fourth of the environment carrying capacity. In such a case,
the involved mechanisms are those of persistence, which require in
particular the state $0$, where no individuals are present, to be
very repulsing  \citep[or equivalently ``unstable",
see][]{sk,bhr1}. It is therefore not surprising to obtain effects
of fragmentation comparable to those described in the existing
literature.

On the other hand, the yield $R(5)$ again reaches a maximum (cf.
Fig. \ref{fig:prop}, b).  As in the quasi-constant-yield case, the
maximum yield is attained for the most fragmented configuration.
Fragmented configurations have a longer perimeter, and thus
provide higher transfer rates into the harvested regions. Indeed,
obtaining higher yields with smaller total populations, for small
values of $E$, implies that the size of the population situated
outside the protected regions is higher for fragmented
configurations; see Note~1. In the quasi-constant case, such a
larger ``unprotected" population would not have implied higher
catches whenever $u$ greater than $\varepsilon$; this explains the
qualitative difference between the two harvesting strategies in
terms of the effects of fragmentation on $P(5)$.

As above discussed, higher values of the effort give a significant
advantage to aggregated configurations, in terms of total population
sizes. This translates into higher yields on aggregated
configurations. Still, these qualitative results do not depend on
the values of diffusion coefficients, in the selected range $[10,100]$, and are not shown.

\subsection{Total population  vs. annual yield}\label{sec_resRP}

\begin{figure}
\centering
\includegraphics*[width=14cm]{colorbar.eps}
\subfigure[]{%
\includegraphics*[width=6.5cm]{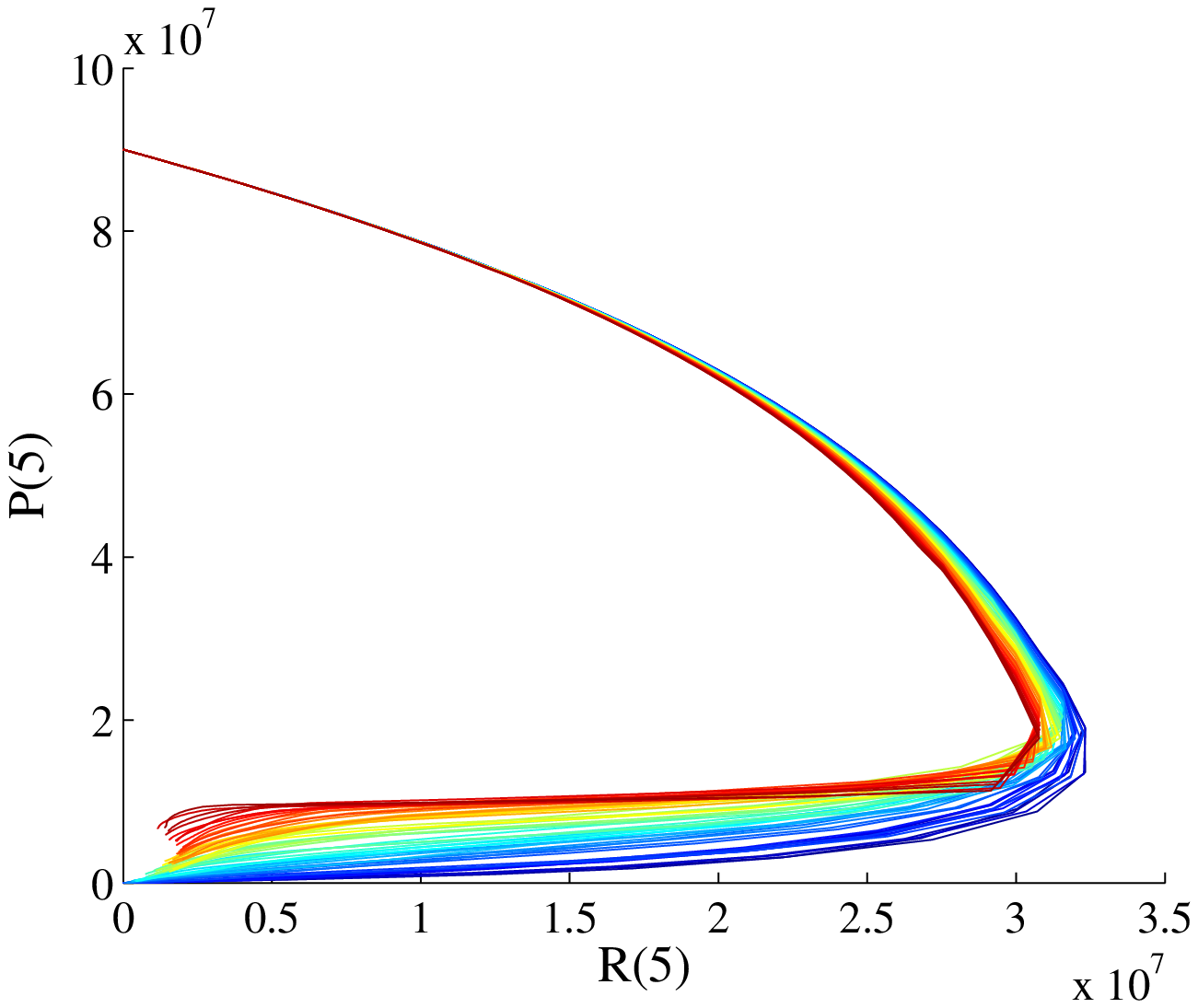}}
\subfigure[]{%
\includegraphics*[width=6.5cm]{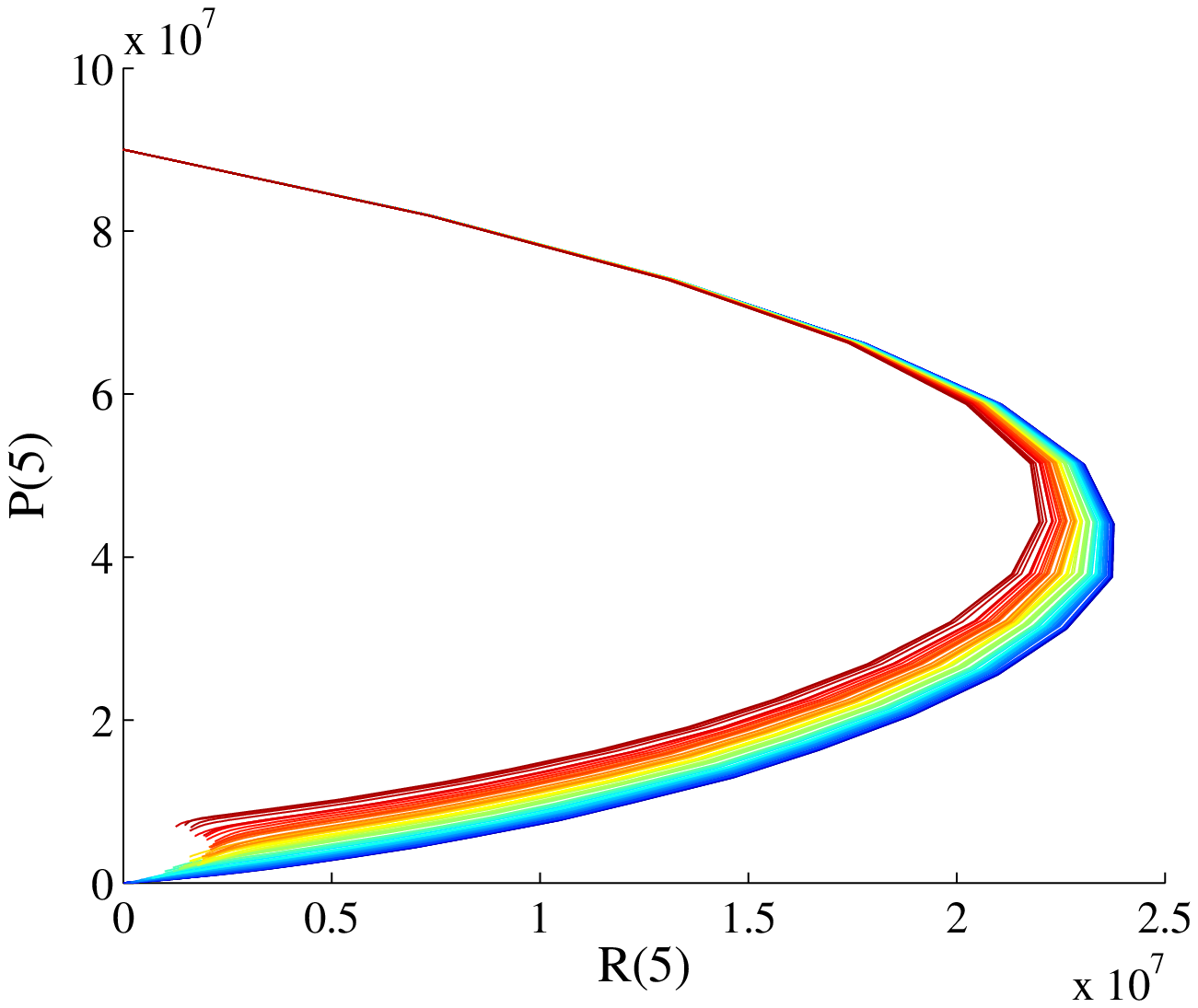}}
\caption{\footnotesize (a): Total population $P(5)$ in $\Omega$
after $5$ years, in terms of the yield during year $4$, $R(5)$, for
the quasi-constant-yield strategy. (b): Total population $P(5)$ vs
yield $R(5)$, for the proportional harvesting strategy. Blue curves
correspond to more fragmented configurations, and
red curves correspond to more aggregated configurations. For these computations, we fixed $D=50$. }%
\label{fig:PR}
\end{figure}

In Fig. \ref{fig:PR},  are depicted the population sizes $P(5)$, in
terms of the yields $R(5)$, and of the level of fragmentation of the
protected region for the quasi-constant-yield and proportional harvesting
strategies (Figs. \ref{fig:PR}, (a) and (b), respectively). The flip
shapes of these diagrams teach us that for a given yield two branch
of disjoint intervals of population size are admissible,
sufficiently far to the left from the ``bending point". The upper
branch corresponds to low harvesting intensities (below that leading
to maximum yield), and the lower branch corresponds to higher
intensities. Remarkably, for each given population size, higher
yields are obtained on more fragmented configurations; the lower the
population size, the higher this effect.

This $\supset$-shape of the $P-R$ diagram supports the idea that it
will be difficult to predict the quantitative effect of
fragmentation in practice, at a fixed yield, without knowledge of
the total population size; a situation which typically arises in
ecological application where such a knowledge is difficult to
achieve.

\subsection{Results for times $t\in (0,20)$}\label{sec_resPt}

\begin{figure}
\centering
\includegraphics*[width=14cm]{colorbar.eps}
\subfigure[]{%
\includegraphics*[width=6.5cm]{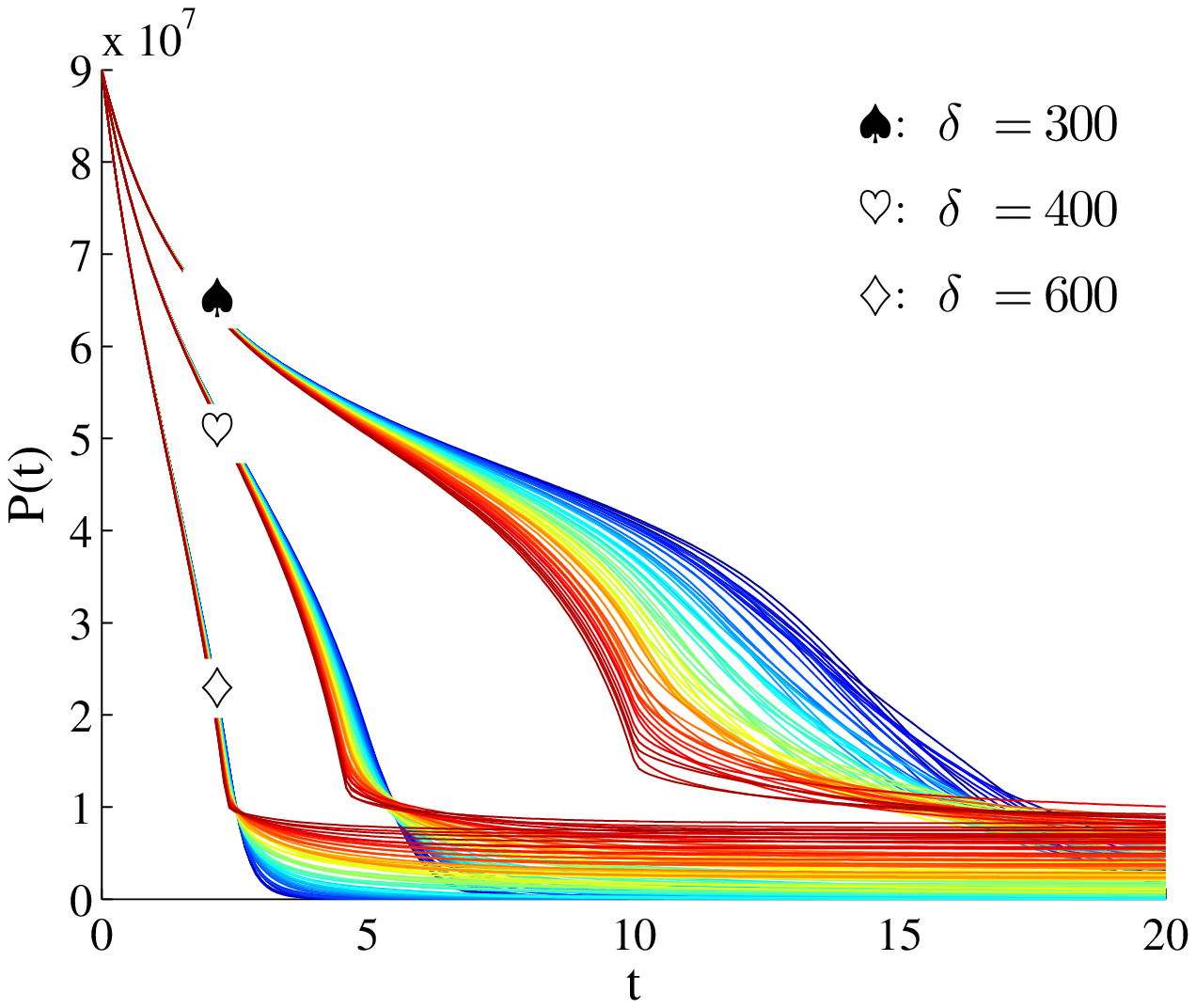}}
\subfigure[]{%
\includegraphics*[width=6.5cm]{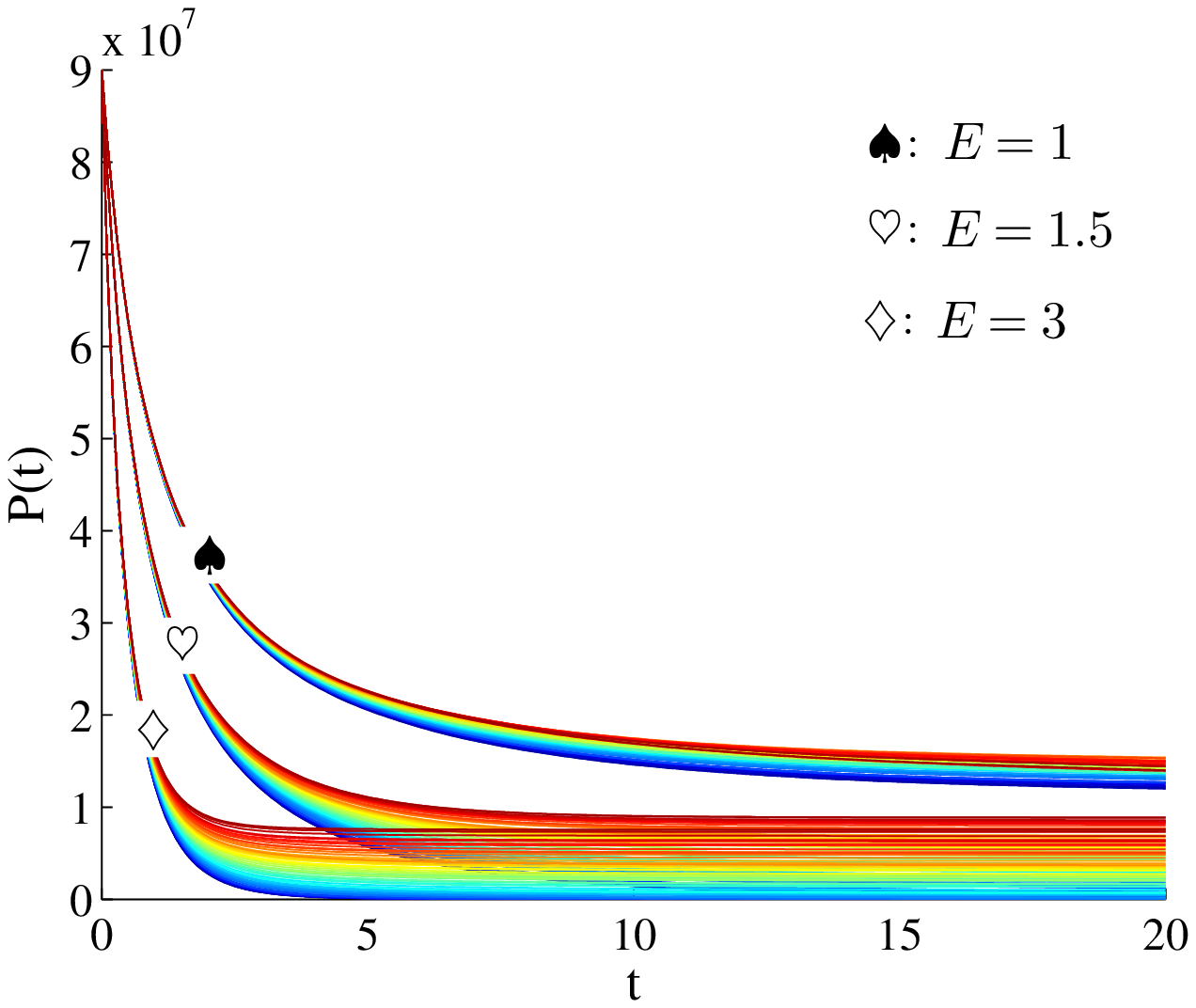}}
\caption{\footnotesize (a): Total population $P(t)$
 in function of time $t$ and in terms of the quota $\delta$, for
the quasi-constant-yield strategy. (b): Total population $P(t)$
 in function of time $t$ and in terms of the effort $E$, for
 the proportional harvesting strategy.  For these computations, we again fixed $D=50$. }%
\label{fig:Pt}
\end{figure}

Fig. \ref{fig:Pt} depicts how the total population $P(t)$ depends on the harvesting intensity and on the aggregation index in function of the time  $t$.

For the quasi-constant-yield harvesting strategy (Fig. \ref{fig:Pt}, a), we observe that the higher the quota, the sooner the inversion of the effects of fragmentation. Thus, the threshold quotas, above which aggregated configurations lead to higher populations, decrease as $t$ increases.

For the proportional harvesting strategy (Fig. \ref{fig:Pt}, b), at each time $t$,  aggregated configurations are still associated with
larger populations. Moreover, the effect of fragmentation tends to increase with time.

In both cases,  the amplitude of the effect of fragmentation depends on the harvesting intensity and on the time point of the analysis. However, those results, and in particular the quasi-constant-yield  case, suggest previous sections results might  still be qualitatively true at times other that $t=5$. This is confirmed by numerical computations (not presented here) which indeed lead to comparable results, but with a shift in the harvesting intensity.

Analytical studies  show that the solution $u(t,x)$ of (\ref{eq_base}) converges to some equilibrium state for both quasi-constant-yield \citep{rc1} and proportional \citep{bhr1} harvesting strategies; however, convergence rates  are not known. Figs. \ref{fig:Pt} (a) and (b)  provide information about these convergence rates: the higher the harvesting intensities, the sooner equilibria are reached. Thus, at  $t=5$,  population sizes may be almost at equilibrium for high harvesting intensities or still farer to reach it for lower intensities.

\section{Concluding remarks}\label{sec_ccl}

In the existing single-species modeling literature, fragmentation
of the favorable region is usually found to be detrimental for
population survival. On the other hand, recent empirical studies
show that fragmentation may have positive or negative effects, and
that positive effects occur more often.

In her review paper on the topic \cite{fahrig2003} proposed two
reasons for negative effects of fragmentation of these favorable
regions. Firstly, favorable patches are too small to sustain a
local population, and secondly, the total perimeter of the
favorable region is large, leading to increased transfer rates
into the unfavorable regions. She also proposed several reasons
for positive effects of fragmentation. Among these reasons,
smaller distance between patches, higher immigration rates into
the patches \citep[see][]{fahrig2004} and positive edge effects
for some species are evoked.

In our work, based on single-species reaction-diffusion models with
harvesting terms, we have captured both positive and negative
effects of fragmentation of the favorable region, interpreted here
as a protected region.

Firstly, for large harvesting terms, we found that aggregated
configurations of the protected region lead to higher population sizes, and
higher yields than fragmented configurations. In that sense, our
results are not contradictory with previous modeling results: under
hostile conditions, the extinction risks are higher, and the chances
of persistence are increased on aggregated configurations. The
reasons evoked by Fahrig for negative  effects of fragmentation then
become paramount in our model at low population sizes.

On the other hand, fragmented protected regions lead to higher
population sizes when a constant number of individuals are removed
per unit of time in the unprotected region
(constant-yield-harvesting). Furthermore, for fixed population
sizes, the two harvesting strategies we studied in this paper lead
to higher yields for fragmented configurations;  cf.  Fig.
\ref{fig:PR}. Such higher yields may stimulate harvesters to slow
down, and may therefore be beneficial for populations.

\cite{sumaila}, through a compartment modeling approach, showed
that transfer rates between favorable and unfavorable regions were
of critical importance for understanding the role of reserves, and
the author emphasized the necessity of a more precise modeling
approach of the transfer rate function. Latter, \cite{tischendorf}
found, via a simulation approach, negative effects of
fragmentation when the probability of individuals to go from the
favorable to the unfavorable region was high, and positive effects
in the opposite case with a high unfavorable to favorable
boundary-crossing probability. In our models, the diffusion
coefficient is spatially-constant, and the boundary-crossing
probabilities are therefore equal in inward and outward
directions. However, net transfer rates from favorable to
unfavorable regions depend on the geometry of the protected
regions, and on the relative population densities inside and
outside the favorable regions. Indeed, individuals moves are
driven by random diffusion, and Green's formula \citep[see
e.g.][]{evans} implies that the instantaneous population flux from
the protected region to the harvested one is:
\par\nobreak\noindent
$$\hbox{Flux}_{\hbox{Protected region}\to \hbox{harvested region}}=-D \int_{\Gamma}\frac{\partial u}{\partial
n} ds,$$where $\Gamma$ denotes the boundary of the protected
region, and $\frac{\partial u}{\partial n}$ is the outward
gradient in population density observed on this boundary. Our work
shows that both the geometry of the protected regions, and the
relative population densities inside and outside the favorable
regions interact to give negative effects of fragmentation for
high harvesting terms,  i.e. when the contrast between protected
and harvested areas is strong, and positive effects for less
contrasted environments.

In multi-species models, fragmentation  per se can alter
interactions among species. As reviewed by \citet{ryall},
predator-prey models predict varying effects of fragmentation on
equilibrium  densities of predator and prey populations, depending
on the specific assumptions of these models. In this paper we have
demonstrated that single-species reaction-diffusion models with
removal terms can support relative effects of fragmentation  per
se. The methods and results could serve as a first step to bridge
the gap between empirical work and modeling, for
reaction-diffusion as well as other single-species models.

\section*{Acknowledgements}

The authors are indebted to an anonymous referee and to Andr\'e Kretzschmar for their valuable comments.
The authors are supported by the French ``Agence Nationale de la Recherche"
within the projects PREFERED, ColonSGS (first author) and URTICLIM (first author).


\begin{thebibliography}{99}

\footnotesize{
\bibitem[Beddington and May(1977)]{bedmay}
 Beddington, J. R., and R. M. May. 1977. Harvesting natural populations in a
randomly fluctuating environment. Science 197:463--465.


\bibitem[Berestycki et al.(2005a)]{bhr1}
Berestycki, H., F. Hamel, and L. Roques. 2005a. Analysis of the periodically
fragmented environment model:\/ {\rm I---}Species persistence. Journal of Mathematical Biology
51:75--113.

\bibitem[Berestycki et al.(2005a)]{bhr2}
Berestycki, H., F. Hamel, and L. Roques. 2005b.  Analysis of the
periodically fragmented environment model:\/ {\rm II---}Biological
invasions and pulsating travelling fronts. Journal de Math\'{e}matiques Pures et Appliqu\'{e}es 84:1101--1146.

\bibitem[Bolker(2003)]{bolker}
Bolker, B. M. 2003. Combining endogenous and exogenous spatial
variability in analytical population models. Theoretical
Population Biology 64:255ñ-270.

\bibitem[Cantrell and Cosner(1989)]{cc0}  Cantrell, R. S., and C. Cosner. 1989. Diffusive logistic equations with
indefinite weights: population models in disrupted environments.
Proceedings of the Royal Society of Edinburgh 112:293--318.

\bibitem[Cantrell and Cosner(2003)]{ccL} Cantrell R. S., and C. Cosner. 2003.  Spatial Ecology via
Reaction-Diffusion Equations. Ser. Math. Comput. Biol., John Wiley
and Sons, Chichester, UK.

\bibitem[Chekroun and Roques(2006)]{cr1}  Chekroun M. D., and L. J. Roques. 2006.
Models of population dynamics under the influence of external
perturbations: mathematical results. Comptes  Rendus Mathematique,
Acadad\'emie des Sciences de Paris, Serie I 343:307--310.

\bibitem[Christensen et al.(2003)]{christ}
Christensen C., S. Gu\'enette, J. J. Heymans, C. J. Walters, R.
Watson, D. Zeller, and D. Pauly. 2003. Hundred-year decline of
North Atlantic predatory fishes. Fish and Fisheries 4:1--24.

\bibitem[El Smaili et al.(2009)]{ehr} El Smaili, M., F. Hamel, and L. Roques. 2009. Homogenization and influence of fragmentation in a biological invasion
model. Discrete and Continuous Dynamical Systems Series A,  25:321--342.

\bibitem[Evans(1998)]{evans} Evans, L. C. 1998. Partial Differential Equations. University of California, Berkeley - AMS.

\bibitem[Fahrig(2003)]{fahrig2003}
Fahrig, L. 2003. Effects of habitat fragmentation on biodiversity.
Annual Review of Ecology, Evolution and Systematics, 34:487--515.


\bibitem[Fisher(1937)]{fi} Fisher, R. A. 1937.  The advance of advantageous genes. Annals of Eugenics, 7:335--369.

\bibitem[Gardner et al.(1987)]{gard}
Gardner, R. H., B. T. Milne, M. G. Turner, and R. V. O'Neill.
1987. Neutral models for the analysis of broad-scale landscape
pattern. Landscape Ecololgy  1:19--28.

\bibitem[Grez et al.(2004)]{fahrig2004} Grez, A., T. Zaviezo, L. Tischendorf, and L. Fahrig. 2004.
A transient, positive effect of habitat fragmentation on insect
population densities. Oecologia 141:444--451

\bibitem[Hale and Verduyn Lunel(1990)]{hale} Hale, J. K., and S. M. Verduyn Lunel. 1990. Averaging in infinite
dimensions. Journal of Integral Equations and Applications
2:463--494.

\bibitem[Holmes(1993)]{Holmes93} Holmes, E.E. 1993. Are Diffusion Models too Simple? A Comparison with Telegraph
Models of Invasion. The American Naturalist 142:779--795.

\bibitem[Huffaker(1958)]{huffaker} Huffaker, C.B. 1958. Experimental Studies on Predation: Dispersion
factors and predator-prey oscillations. Hilgardia  27:343--383.


\bibitem[Keitt(2000)]{keitt} Keitt, T. H. 2000. Spectral representation of neutral landscapes.
Landscape Ecology  15:479--494.

\bibitem[Kinezaki et al.(2003)]{kkts}  Kinezaki, N., K.  Kawasaki, F.  Takasu, and N. Shigesada. 2003.
Modeling biological invasion into periodically fragmented
environments. Theoretical Population Biology 64:291--302.

\bibitem[Kolmogorov et al.(1937)]{kpp}  Kolmogorov, A. N., I. G. Petrovsky, and N. S. Piskunov. 1937.
Etude de l'\'equation de la diffusion avec croissance de la
quantit\'e de mati\`ere et son application \`a un probl\`eme
biologique, Bulletin Universit\'e d'\'Etat \`a Moscou (Bjul.
Moskowskogo Gos. Univ.), S\'erie internationale A  1:1--26.


\bibitem[Mandelbrot(1982)]{mandel} Mandelbrot, B. B. 1982.  The Fractal Geometry of Nature, W. H. Freeman,
New York.

\bibitem[Murray(2002)]{mu} Murray, J. D. 2002. Mathematical biology. Third edition, Interdisciplinary Applied Mathematics 17, Springer-Verlag, New York.

\bibitem[Neubert(2003)]{neub} Neubert, M. G. 2003. Marine reserves and optimal harvesting, Ecology Letters 6:843--849.

\bibitem[Okubo and Levin(2002)]{oku}
 Okubo, A., and S. A. Levin. 2002.  Diffusion and Ecological
Problems---Modern Perspectives. 2nd ed., Springer-Verlag, New
York.

\bibitem[Oruganti et al.(2002)]{shi2}
Oruganti, S.,  R. Shivaji, and J. Shi. 2002.  Diffusive logistic
equation with constant effort harvesting,\/ {\rm I:} Steady
states, Transactions of the American Mathematical Society
354:3601--3619.

\bibitem[Pearson and Blakeman(1906)]{pm}  Pearson, K., and J. Blakeman. 1906.  Mathematical Contributions of the Theory of Evolution, XV, A Mathematical Theory of Random Migration. Drapersi Company Research Mem. Biometrics Series~III, Dept. Appl. Meth. Univ. College, London.

\bibitem[Robinson and Bodmer(1999)]{rob2}
 Robinson J. G., and R. E. Bodmer. 1999. Towards wildlife management in
tropical forests. Journal of Wildlife Management  63:1--13.

\bibitem[Robinson and Redford(1991)]{rob1}
 Robinson J. G. and K. H. Redford. 1991.  Sustainable harvest of
neo-tropical mammals. Pages 415--429 in J. G. Robinson and K. H.
Redford, eds. Neo-Tropical Wildlife Use and Conservation, Chicago
University Press, Chicago, IL.

\bibitem[Roques and Chekroun(2007)]{rc1}  Roques, L., and  M. D. Chekroun. 2007. On population resilience to external
perturbations. SIAM Journal on Applied Mathematics 68:133--153.

\bibitem[Roques and Hamel(2007)]{rh1}
 Roques, L., and F. Hamel. 2007.  Mathematical analysis of the optimal
habitat configurations for species persistence. Mathematical
Biosciences  210:34--59.

\bibitem[Roques and Stoica(2007)]{rs1}
Roques, L., and R. Stoica. 2007.  Species persistence decreases
with habitat fragmentation: An analysis in periodic stochastic
environments. Journal of Mathematical Biology 55:189--205.

\bibitem[Ryall and Fahrig(2006)]{ryall}  Ryall, K. L., and L. Fahrig. 2006. Response of predators to loss and fragmentation of prey habitat: a review of theory. Ecology  87:1086--1093.

\bibitem[Saunders et al.(1991)]{shm} Saunders, D.A., R. J. Hobbs, and C. R. Margules. 1991. Biological
consequences of ecosystem fragmentation: a review. Conservation
Biology 51:18--32.

\bibitem[Shigesada and Kawasaki(1997)]{sk}
N. Shigesada and K. Kawasaki,. 1997.  Biological Invasions: Theory
and Practice. Oxford Series in Ecology and Evolution, Oxford
University Press, Oxford, UK.

\bibitem[Skellam(1951)]{ske} Skellam, J. G. 1951. Random dispersal in theoretical populations. Biometrika 38:196-218.

\bibitem[Stephens et al.(2002)]{steph}
Stephens, P. A., F. Frey-Roos, W. Arnold, and W. J. Sutherland.
2002. Sustainable exploitation of social species: A test and
comparison of models. Journal of Applied Ecology 39:629--642.

\bibitem[Sumaila(1998)]{sumaila} Sumaila, U. R. 1998. Protected marine reserves as fisheries management tools: a bioeconomic
analysis, Fisheries Research 37:287--296.



\bibitem[Tischendorf et al.(2005)]{tischendorf} Tischendorf, L.,  A. Grez, T. Zaviezo, and L. Fahrig. 2005.
Mechanisms affecting population density in fragmented habitat.
Ecology and Society 10:7.


\bibitem[Turchin(1998)]{tur}
Turchin, P. 1998.  Quantitative Analysis of Movement: Measuring
and Modeling Population Redistribution in Animals and Plants,
Sinauer Associates, Sunderland, MA.

}

\end{thebibliography}
\end{document}